# Robust Predictor Feedback for Discrete-Time Systems with Input Delays

**Iasson Karafyllis*** **and Miroslav Krstic****

*Dept. of Environmental Eng., Technical University of Crete,
73100, Chania, Greece, email: ikarafyl@enveng.tuc.gr

**Dept. of Mechanical and Aerospace Eng., University of California,
San Diego, La Jolla, CA 92093-0411, U.S.A., email: krstic@ucsd.edu

## Abstract

This work studies the design problem of feedback stabilizers for discrete-time systems with input delays. A backstepping procedure is proposed for disturbance-free discrete-time systems. The feedback law designed by using backstepping coincides with the predictor-based feedback law used in continuous-time systems with input delays. However, simple examples demonstrate that the sensitivity of the closed-loop system with respect to modeling errors increases as the value of the delay increases. The paper proposes a Lyapunov redesign procedure which can minimize the effect of the uncertainty. Specific results are provided for linear single-input discrete-time systems with multiplicative uncertainty. The feedback law that guarantees robust global exponential stability is a nonlinear, homogeneous of degree 1 feedback law.

## 1. Introduction

Continuous-time systems with input delays have been studied extensively in the literature (see [1,17,18] and the references therein). However, discrete-time systems with input delays have been rarely studied. In many aspects the results in discrete-time systems are complementary to the results obtained for continuous-time systems. The papers [9,23] tried to extend the Smith predictor design for nonlinear discrete-time systems, which are feedback linearizable. The works [4,5,6,19,21] have provided results for discrete-time linear systems with input or measurement delays. In most cases where the delay is time-varying, the delay is considered to be an unknown perturbation of a nominal value and a predictor-based design is implemented.

It should be noted that discrete-time systems with input delays is a special class of systems with characteristics that are not met in continuous-time systems. The following list summarizes some of the differences between discrete-time systems and continuous-time systems:

1) Discrete-time systems with input delays of the form

$$x(t+1) = f(x(t), u(t-r))$$
$$x(t) \in \Re^n, u(t) \in \Re^m \tag{1.1}$$

where $f : \Re^n \times \Re^m \to \Re^n$ is a mapping with $f(0,0) = 0$ and $r \geq 1$ is an integer, or uncertain discrete-time systems with input delays of the form

$$x(t+1) = F(d(t), x(t), u(t-r))$$
$$x(t) \in \Re^n, u(t) \in \Re^m, d(t) \in D \tag{1.2}$$

where $D \subseteq \Re^l$ is a non-empty set, $d \in D$ denotes the vector of unknown time-varying paramenters (uncertainties), $F : D \times \Re^n \times \Re^m \to \Re^n$ is a mapping with $F(d,0,0) = 0$ for all $d \in D$ and $r \geq 1$ is an integer, are always forward complete. This is a major difference with continuous-time systems, where forward completeness is not guaranteed.

2) The closed-loop system (1.1) with the feedback law

$$u(t) = K(x(t), u(t-r), ..., u(t-1)) \tag{1.3}$$

where $K : \Re^n \times \Re^{rm} \to \Re^m$ is a mapping with $K(0,0,...,0) = 0$, or the closed-loop system (1.2) with the feedback law (1.3), are again finite-dimensional discrete-time systems, i.e., the feedback law (1.3) preserves the qualitative characteristics of the system. This is in sharp contrast with continuous-time systems: as shown in the literature (see [1,17,18] and the references therein) the closed-loop system of a finite-dimensional control system with input delays and a feedback law that depends on the history of the input is a system with very different characteristics (it is not a finite-dimensional system). To see why the character of the system is preserved for discrete-time systems, we should note that system (1.1) is equivalent to the system:

$$x(t+1) = f(x(t), y_1(t))$$
$$y_i(t+1) = y_{i+1}(t) \quad , \quad i = 1,...,r$$
$$x(t) \in \Re^n, y_i(t) \in \Re^m \ (i = 1,...,r), y_{r+1}(t) = u(t) \in \Re^m \tag{1.4}$$

and system (1.2) is equivalent to the system:

$$x(t+1) = F(d(t), x(t), y_1(t))$$
$$y_i(t+1) = y_{i+1}(t) \quad , \quad i = 1,...,r$$
$$x(t) \in \Re^n, y_i(t) \in \Re^m \ (i = 1,...,r),$$
$$y_{r+1}(t) = u(t) \in \Re^m, d(t) \in D \tag{1.5}$$

Indeed, it is straightforward to verify that the following equalities hold:

$$y_i(t) = u(t-r-1+i), \text{ for } i = 1,...,r \text{ and } t \geq r+1-i \tag{1.6}$$

Therefore, the closed-loop system (1.1) with the feedback law (1.3) is equivalent to the closed-loop system (1.4) with

$$u(t) = K(x(t), y_1(t), ..., y_r(t)) \tag{1.7}$$

Clearly, the closed-loop system (1.4) with (1.7) is a finite-dimensional discrete-time system. Similarly, the closed-loop system (1.2) with the feedback law (1.3) is equivalent to the closed-loop system (1.5) with (1.7), which again is an uncertain finite-dimensional discrete-time system.

3) The equivalent description of the discrete-time systems (1.1) or (1.2), i.e., systems (1.4) and (1.5), respectively, have a specific structure: they are composed from a nonlinear component and a cascade of "sum-ators". This specific structure can be exploited in order to design the feedback law (1.7) efficiently by means of backstepping: this feature is absent in the analysis of continuous-time systems. The backstepping feedback design for discrete-time systems was studied in [11,12,22] and is in the same spirit of the backstepping feedback design for continuous-time systems (see [16]).

4) The following implication holds: if the feedback law (1.3) stabilizes the equilibrium point for (1.1) or (1.2) then the feedback law

$$u(t) = K(x(t-r), u(t-r), ..., u(t-1)) \tag{1.8}$$

stabilizes the equilibrium point for $x(t+1) = f(x(t), u(t))$ or $x(t+1) = F(d(t), x(t), u(t))$, respectively. In other words, input delays and measurement delays can be treated in the same way: this convenient feature is not completely true in continuous-time systems (see the discussion in [13]).

This work is devoted to the answer of the following question: how can we design a feedback law of the form (1.7) so that $0 \in \Re^n \times \Re^{rm}$ is (robustly) globally asymptotically stable for the corresponding closed-loop system. Notice that we consider nonlinear systems with uncertainties: system (1.2) is a discrete-time system with vanishing perturbations (multiplicative uncertainties). The structure of this paper is as follows:

- Section 2 describes a backstepping solution which exploits the specific structure of systems (1.4) and (1.5). It is shown that the backstepping solution is a direct extension of the predictor-based approach, which has already been described for continuous-time systems. Our approach does not require the assumption of open-loop stability, which is present in [9,23].

- As expected, simple examples show that the sensitivity of the closed-loop system with respect to uncertainties is magnified as the value of the delay $r \geq 1$ increases. For this reason, it is important to redesign the stabilizing feedback after applying the backstepping approach in order to reduce the sensitivity: this is the topic of Section 3. Specific results for the Lyapunov redesign are provided for linear single-input systems with vanishing perturbations (Theorem 3.2). Explicit inequalities allow the determination of the lowest upper bound for the magnitude of the uncertainty for which robust global exponential stability holds for the closed-loop system. A simple example shows the importance of the Lyapunov redesign procedure (Example 3.4).

It should be noted that Lyapunov redesign is a well-known procedure for nonlinear continuous-time systems (see [15]). Recently, Lyapunov redesign has been used extensively in sampled-data feedback design (see [7,8,20]).

*Notation.* Throughout the paper we adopt the following notation:

* For a vector $x \in \Re^n$ we denote by $|x|$ its usual Euclidean norm, by $x'$ its transpose. For a real matrix $A \in \Re^{n\times m}$, $A' \in \Re^{m\times n}$ denotes its transpose and $|A| := \sup\{\,|Ax|\,; x \in \Re^n\,, |x| = 1\,\}$ is its induced norm. $I \in \Re^{n\times n}$ denotes the identity matrix.
* $\Re_+$ denotes the set of non-negative real numbers.

Throughout the paper, we assume that the mapping $f : \Re^n \times \Re^m \to \Re^n$ appearing in the right hand side of (1.1) is continuous. Moreover, we assume that the mapping $F : D \times \Re^n \times \Re^m \to \Re^n$ appearing in the right hand side of (1.2) is continuous. The notions of (robust) global asymptotic stability and robust global exponential stability employed in this work are the standard notions described in [10,12,14].

## 2. A Backstepping Solution

We assume, as in the continuous-time case, that system (1.1) with $r = 0$ (the delay-free version of (1.1)) is stabilizable, i.e., we make the following assumption.

**(H1)** *There exists a continuous function* $k : \Re^n \to \Re^m$ *with* $k(0) = 0$ *such that* $0 \in \Re^n$ *is Globally Asymptotically Stable (GAS) for the closed-loop system (1.1) with* $r = 0$ *and* $u(t) = k(x(t))$.

In order to be able to address the stabilization problem for (1.1) with $r > 0$, or equivalently the stabilization problem for (1.4), we need the following technical lemma.

**Lemma 2.1 (the backstepping lemma for discrete-time systems):** *Suppose that (H1) holds. Let a continuous, positive definite and radially unbounded function* $V : \Re^n \to \Re_+$ *and a constant* $\lambda \in [0,1)$ *be such that the following inequality holds:*

$$V(f(x, k(x))) \le \lambda V(x)\text{, for all } x \in \Re^n \tag{2.1}$$

*Then the following hold:*

**(i)** $0 \in \Re^n \times \Re^m$ *is GAS for the closed-loop (1.4) with* $r = 1$ *and* $u(t) = k\big(f(x(t), y_1(t))\big)$,

**(ii)** *for every* $a \in K_\infty$ *and for every constant* $c > \frac{1}{1-\lambda}$, *the function* $\bar{V} : \Re^n \times \Re^m \to \Re_+$ *defined by*

$$\bar{V}(x, y_1) := V(x) + cV(f(x, y_1)) + a\big(|y_1 - k(x)|\big) \tag{2.2}$$

*is a continuous, positive definite and radially unbounded function that satisfies:*

$$\bar{V}(f(x, y_1), k(f(x, y_1))) \le \left(\lambda + c^{-1}\right)\bar{V}(x, y_1)\text{, for all } (x, y_1) \in \Re^n \times \Re^m \tag{2.3}$$

**Proof:** (i) is a direct consequence of (ii), the fact that $\lambda + c^{-1} < 1$ (which is a consequence of $c > \frac{1}{1-\lambda}$) and the Lyapunov theorem for discrete-time systems (see [10,12,14]). Therefore, we focus on proving (ii).

Continuity of $\bar{V} : \Re^n \times \Re^m \to \Re_+$ is a direct consequence of continuity of $V : \Re^n \to \Re_+$, $f : \Re^n \times \Re^m \to \Re^n$ and $a \in K_\infty$. The fact that $\bar{V} : \Re^n \times \Re^m \to \Re_+$ as defined by (2.2) is positive definite is a direct consequence of definition (2.2) and the fact that $k(0) = 0$. In order to show that $\bar{V} : \Re^n \times \Re^m \to \Re_+$ is radially unbounded, it suffices to show that the set

$$S_M := \left\{ (x, y_1) \in \Re^n \times \Re^m : \bar{V}(x, y_1) \le M \right\} \tag{2.4}$$

is bounded for every $M \ge 0$. Indeed, by virtue of definition (2.2), we conclude that for each $(x, y_1) \in S_M$ it holds that $V(x) \le M$. Since $V : \Re^n \to \Re_+$ is radially unbounded, it follows that the component $x \in \Re^n$ of the vector $(x, y_1) \in S_M$ is bounded, i.e., there exists $R \ge 0$ such that $|x| \le R$ for all $(x, y_1) \in S_M$. Since $k : \Re^n \to \Re^m$ is continuous, it follows that there exists $\tilde{R} \ge 0$ such that $|k(x)| \le \tilde{R}$, for all $(x, y_1) \in S_M$. Finally, notice that definitions (2.2), (2.4) allow us to conclude the following inequality for all $(x, y_1) \in S_M$:

$$|y_1 - k(x)| \le a^{-1}(M) \tag{2.5}$$

Combining (2.5) with the inequality $|k(x)| \le \tilde{R}$, for all $(x, y_1) \in S_M$, we obtain $|y_1| \le a^{-1}(M) + \tilde{R}$, for all $(x, y_1) \in S_M$. In other words, the component $y_1 \in \Re^m$ of the vector $(x, y_1) \in S_M$ is bounded. Therefore, $S_M$ as defined by (2.4) is bounded for every $M \ge 0$.

We are left with the task of proving (2.3). Definition (2.2) implies:

$$\bar{V}(f(x, y_1), k(f(x, y_1))) = V(f(x, y_1)) + cV(f(f(x, y_1), k(f(x, y_1)))), \text{ for all } (x, y_1) \in \Re^n \times \Re^m \tag{2.6}$$

Using (2.1) with $f(x, y_1)$ in place of $x \in \Re^n$ and (2.6), we obtain:

$$\bar{V}(f(x, y_1), k(f(x, y_1))) \le (1 + c\lambda) V(f(x, y_1)), \text{ for all } (x, y_1) \in \Re^n \times \Re^m \tag{2.7}$$

Finally, notice that definition (2.2) implies $V(f(x, y_1)) \le \frac{1}{c} \bar{V}(x, y_1)$ for all $(x, y_1) \in \Re^n \times \Re^m$. The previous inequality in conjunction with (2.7) gives (2.3).

The proof is complete. ◁

The reader should notice at this point that the existence of a continuous, positive definite and radially unbounded function $V : \Re^n \to \Re_+$ and a constant $\lambda \in [0,1)$ satisfying (2.1) is a direct consequence of assumption (H1) and Proposition 3.1 in [12].

Applying Lemma 2.1 inductively, allows us to construct a globally stabilizing feedback for system (1.4). We define the extended vectors

$$z_0 = x, \; z_i = (x, y_1, ..., y_i), \text{ for } i = 1, ..., r+1 \tag{2.8}$$

and the vector fields $F_i(z_i)$, $i = 0,1,...,r+1$ by the recursive formula:

$$F_{i+1}(z_{i+1}) := f\big(F_i(z_i), y_{i+1}\big) \text{ for } i = 0,...,r \text{ and } F_0(x) := x \tag{2.9}$$

We notice that the following formulae hold:

$$F_i(f(x, y_1), y_2, ..., y_{i+1}) = F_{i+1}(x, y_1, ..., y_{i+1}) = F_{i+1}(z_{i+1}) \text{ for } i = 1,...,r \tag{2.10}$$

We are now ready to state and prove a result which deals with the global stabilization of (1.4). The proof of Lemma 2.2 can be made by induction. However, for clarity purposes, we will provide a different (and more direct) proof for Lemma 2.2, which uses the identities (2.9) and (2.10).

**Lemma 2.2:** *Consider system (1.4) with $r \geq 2$ and suppose that (H1) holds. Let a continuous, positive definite and radially unbounded function $V : \Re^n \to \Re_+$ and a constant $\lambda \in [0,1)$ be such that (2.1) holds. Then the following hold:*
**(i)** $0 \in \Re^n \times \Re^m$ *is GAS for the closed-loop (1.4) and $u(t) = k\big(F_r(z_r(t))\big)$,*
**(ii)** *for every $a_i \in K_\infty$ ($i = 1,...,r$) with $a_i(s) \leq a_{i+1}(s)$ for all $i = 1,...,r-1$ and $s \geq 0$, for every constant $c > \frac{1}{1-\lambda}$, the function $\bar{V} : \Re^n \times \Re^{rm} \to \Re_+$ defined by*

$$\bar{V}(z_r) := \sum_{i=0}^{r} c^i V\big(F_i(z_i)\big) + \sum_{i=1}^{r} c^i a_i\Big(\big|y_i - k\big(F_{i-1}(z_{i-1})\big)\big|\Big) \tag{2.11}$$

*is a continuous, positive definite and radially unbounded function that satisfies:*

$$\bar{V}(f(x, y_1), y_2, ..., y_r, k(F_r(z_r))) \leq \big(\lambda + c^{-1}\big)\bar{V}(z_r), \text{ for all } z_r = (x, y_1, ..., y_r) \in \Re^n \times \Re^{rm} \tag{2.12}$$

**Proof:** (i) is a direct consequence of (ii), the fact that $\lambda + c^{-1} < 1$ (which is a direct consequence of the fact that $c > \frac{1}{1-\lambda}$) and Proposition 2.3 in [12]. The fact that the function $\bar{V} : \Re^n \times \Re^{rm} \to \Re_+$ is a continuous, positive definite and radially unbounded function can be proved in exactly the same way as in the proof of Lemma 2.1. We show next the validity of inequality (2.12).

Using definition (2.11), we obtain:

$$\begin{aligned} \bar{V}(f(x, y_1), y_2, ..., y_r, u) &= c^r V\big(F_r(f(x, y_1), y_2, ..., y_r, u)\big) + \sum_{i=0}^{r-1} c^i V\big(F_i(f(x, y_1), y_2, ..., y_{i+1})\big) \\ &+ \sum_{i=1}^{r} c^i a_i\Big(\big|y_{i+1} - k\big(F_{i-1}(f(x, y_1), y_2, ..., y_i)\big)\big|\Big) \end{aligned} \tag{2.13}$$

The identities (2.10) in conjunction with (2.13) give:

$$\begin{aligned} \bar{V}(f(x, y_1), y_2, ..., y_r, u) &= c^r V\big(F_{r+1}(z_{r+1})\big) + \frac{1}{c}\sum_{i=1}^{r} c^i V\big(F_i(z_i)\big) \\ &+ \sum_{i=1}^{r-1} c^i a_i\Big(\big|y_{i+1} - k\big(F_i(z_i)\big)\big|\Big) + c^r a_r\Big(\big|u - k\big(F_r(z_r)\big)\big|\Big) \end{aligned} \tag{2.14}$$

Using identity (2.9) for $i = r$ and inequality (2.1) with $F_r(z_r)$ in place of $x \in \Re^n$, in conjunction with (2.14), we obtain for $u = k(F_r(z_r))$:

$$\begin{aligned} \bar{V}(f(x, y_1), y_2, ..., y_r, u) &= c^r V(f(F_r(z_r), u)) + \frac{1}{c}\sum_{i=1}^{r} c^i V(F_i(z_i)) \\ &+ c^r a_r(|u - k(F_r(z_r))|) + \sum_{i=1}^{r-1} c^i a_i(|y_{i+1} - k(F_i(z_i))|) \\ &\leq c^r \lambda V(F_r(z_r)) + \frac{1}{c}\sum_{i=1}^{r} c^i V(F_i(z_i)) + \sum_{i=1}^{r-1} c^i a_i(|y_{i+1} - k(F_i(z_i))|) \end{aligned} \tag{2.15}$$

Inequality (2.15) in conjunction with the equality $\frac{1}{c}\sum_{i=1}^{r} c^i V(F_i(z_i)) = \frac{1}{c}\bar{V}(z_r) - \frac{1}{c}V(x) - \frac{1}{c}\sum_{i=1}^{r} c^i a_i(|y_i - k(F_{i-1}(z_{i-1}))|)$ (which is a direct consequence of definition (2.11)) give:

$$\begin{aligned} \bar{V}(f(x, y_1), y_2, ..., y_r, k(F_r(z_r))) &\leq c^r \lambda V(F_r(z_r)) + \frac{1}{c}\bar{V}(z_r) - \frac{1}{c}V(x) \\ &- \sum_{i=1}^{r} c^{i-1} a_i(|y_i - k(F_{i-1}(z_{i-1}))|) + \sum_{i=2}^{r} c^{i-1} a_{i-1}(|y_i - k(F_{i-1}(z_{i-1}))|) \end{aligned} \tag{2.16}$$

Finally, inequality (2.12) is a consequence of (2.16), the fact that $a_i(s) \leq a_{i+1}(s)$ for all $i = 1, ..., r-1$, $s \geq 0$ and the inequality $c^r V(F_r(z_r)) \leq \bar{V}(z_r)$ (which is a direct consequence of definition (2.11)). The proof is complete. ◁

The globally stabilizing feedback laws that are proposed by Lemma 2.1 (for the case $r = 1$) or Lemma 2.2 (for the case $r \geq 2$) are feedback laws which are based on prediction schemes. Indeed, we may verify that the implementation of the feedback law $u(t) = k(f(x(t), y_1(t)))$ for system (1.4) with $r = 1$ guarantees $u(t) = k(x(t+1))$ for all $t \geq 0$ and that the implementation of the feedback law $u(t) = k(F_r(z_r(t)))$ for system (1.4) with $r \geq 2$ guarantees $u(t) = k(x(t+r))$ for all $t \geq 0$. Therefore, the prediction-based control schemes described in [1,13,17,18] for continuous-time systems are modified in an obvious way for discrete-time systems.

For the case of continuous-time systems, the predictor-based control scheme suffers from sensitivity with respect to modeling errors: the sensitivity with respect to modeling errors tends to be larger as the value of the input delay increases. This feature is present in discrete-time systems as well. The following example shows how the sensitivity with respect to vanishing perturbations is magnified as the value of the input delay becomes higher and higher.

**Example 2.3:** Consider the scalar discrete-time system

$$\begin{aligned} &x(t+1) = x(t) + d(t)x(t) + u(t-r) \\ &x(t) \in \Re, d(t) \in D = [-a, a], u(t) \in \Re \end{aligned} \tag{2.17}$$

where $a \geq 0$ and $r \geq 0$ is an integer. System (2.17) is equivalent to the system

$$
\begin{aligned}
&x(t+1)=x(t)+d(t)x(t)+y_1(t)\\
&y_i(t+1)=y_{i+1}(t) \quad , \quad i=1,...,r\\
&x(t)\in\Re^n \, , y_i(t)\in\Re\left(i=1,...,r\right),\\
&y_{r+1}(t)=u(t)\in\Re , d(t)\in[-a,a]
\end{aligned}
\tag{2.18}
$$

The nominal value for the uncertainty $d\in[-a,a]$ is $d\equiv 0$. Based on the nominal value for the uncertain parameter, we obtain the disturbance-free discrete-time system:

$$
\begin{aligned}
&x(t+1)=x(t)+y_1(t)\\
&y_i(t+1)=y_{i+1}(t) \quad , \quad i=1,...,r\\
&x(t)\in\Re^n \, , y_i(t)\in\Re\left(i=1,...,r\right),\\
&y_{r+1}(t)=u(t)\in\Re , d(t)\in[-a,a]
\end{aligned}
\tag{2.19}
$$

The assumptions of Lemma 2.1 (for the case $r=1$) and the assumptions of Lemma 2.2 (for the case $r\geq 2$) hold for system (2.19) with $k(x)=-x$, $V(x)=x^2$ and $\lambda=0$. Therefore, a globally stabilizing feedback for (2.19) is

$$
u(t)=-\left(x(t)+\sum_{i=1}^{r}y_i(t)\right) \tag{2.20}
$$

If the feedback law (2.20) is applied to the uncertain system (2.18), then the following question arises:

"For what values of $a\geq 0$, is $0\in\Re\times\Re^r$ robustly globally asymptotically stable for the closed-loop system (2.18) with (2.20)?"

First, we notice that the closed-loop system (2.18) with (2.20) has a constant and non-zero solution when the constant sequence $d(t)\equiv\frac{1}{r+1}$ is applied and the initial condition satisfies $y_1(0)=...=y_r(0)=-\frac{x(0)}{r+1}\neq 0$. Therefore, it is clear that a necessary condition for robust global asymptotic stability of the equilibrium point $0\in\Re\times\Re^r$ for the closed-loop system (2.18) with (2.20) is

$$
a<\frac{1}{r+1} \tag{2.21}
$$

In order to obtain a sufficient condition for robust global asymptotic stability, we use the knowledge of a Lyapunov function for the closed-loop system (2.18) with (2.20) with $d\equiv 0$. Notice that the function $F_i(z_i)$ defined by (2.9) is given by $F_i(z_i)=x+y_1+...+y_i$ for $i=1,...,r+1$. Therefore, Lemma 2.1 or Lemma 2.2 guarantees that the function

$$
\overline{V}(x,y_1,...,y_r):=x^2+(1+\varphi)\sum_{i=1}^{r}c^i\left(x+y_1+...+y_i\right)^2 \tag{2.22}
$$

with $c>1$, $\varphi>0$, is a Lyapunov function for the closed-loop system (2.18) with (2.20) with $d\equiv 0$. Notice that the function defined by (2.22) corresponds to formula (2.11) with $a_i(s)\equiv\varphi s^2$ for $i=1,...,r$, $k(x)=-x$, $V(x)=x^2$ and $\lambda=0$. However, it is clear that the function defined by (2.22) is a continuous, positive definite and radially unbounded function for all $c>0$ and $\varphi>-1$. Therefore, we consider the function defined by (2.22) for $c>0$ and $\varphi>-1$.

For $r=1$, we get for $u=-x-y_1$ and for all $(x,y_1)\in\Re^2$, $d\in\Re$:

$$\begin{aligned}\bar{V}(x+dx+y_1,u) &= (x+dx+y_1)^2+(1+\varphi)c(x+dx+y_1+u)^2 \\ &= (x+y_1)^2+\left(1+(1+\varphi)c\right)d^2x^2+2dx(x+y_1)\end{aligned} \tag{2.23}$$

Completing the squares, we obtain from (2.23) for all $\varepsilon>0$, $(x,y_1)\in\Re^2$ and $|d|\le a$:

$$\bar{V}(x+dx+y_1,u)\le(1+\varepsilon^{-1})(x+y_1)^2+\left(1+\varepsilon+(1+\varphi)c\right)a^2x^2 \tag{2.24}$$

It follows from (2.24) and (2.22) with $r=1$, that there exists $\sigma\in[0,1)$ such that $\bar{V}(x+dx+y_1,u)\le\sigma\bar{V}(x,y_1)$ holds for all $(x,y_1)\in\Re^2$ and $|d|\le a$, provided that there exists $\varepsilon>0$ so that the following inequalities hold:

$$\begin{aligned}&\left(1+\varepsilon+(1+\varphi)c\right)a^2<1 \\ &1+\varepsilon^{-1}<(1+\varphi)c\end{aligned} \tag{2.25}$$

Inequalities (2.25) hold for certain $\varepsilon>0$, provided that the following inequality holds:

$$a^2<\frac{c(1+\varphi)-1}{c^2(1+\varphi)^2} \tag{2.26}$$

The greatest value for the fraction $\frac{c(1+\varphi)-1}{c^2(1+\varphi)^2}$ is obtained for $c(1+\varphi)=2$. Therefore, using Proposition 2.3 in [12], we can conclude that the equilibrium point $0\in\Re\times\Re$ for the closed-loop system (2.18) with (2.20) and $r=1$ is robustly globally asymptotically stable, provided that $a<1/2$.

Next, we consider the case $r\ge 2$. We get for $u=-x-y_1-\ldots-y_r$ and for all $(x,y_1,\ldots,y_r)\in\Re^{r+1}$, $d\in\Re$:

$$\begin{aligned}\bar{V}(x+dx+y_1,y_2,\ldots,y_r,u) &= (x+dx+y_1)^2+(1+\varphi)\sum_{i=1}^{r}c^i\left(x+dx+y_1+\ldots+y_i+y_{i+1}\right)^2 \\ &= \left(1+(1+\varphi)\sum_{i=1}^{r}c^i\right)d^2x^2+(x+y_1)^2+(1+\varphi)\sum_{i=1}^{r-1}c^i\left(x+y_1+\ldots+y_i+y_{i+1}\right)^2 \\ &\quad+2dx(x+y_1)+2(1+\varphi)dx\sum_{i=1}^{r-1}c^i\left(x+y_1+\ldots+y_i+y_{i+1}\right)\end{aligned} \tag{2.27}$$

Using the inequality $2dx(x+y_1)\le\varepsilon_1 d^2x^2+\varepsilon_1^{-1}(x+y_1)^2$, which holds for every $\varepsilon_1>0$, $(x,y_1)\in\Re^2$, $d\in\Re$ and the inequalities $2dx\left(x+y_1+\ldots+y_i+y_{i+1}\right)\le\varepsilon_2 d^2x^2+\varepsilon_2^{-1}\left(x+y_1+\ldots+y_i+y_{i+1}\right)^2$, which hold for all $i=1,\ldots,r-1$, $\varepsilon_2>0$, $(x,y_1,\ldots,y_r)\in\Re^{r+1}$, $d\in\Re$, we obtain from (2.27) for all $|d|\le a$:

$$\begin{aligned}\bar{V}(x+dx+y_1,y_2,\ldots,y_r,u) &\le\left(1+\varepsilon_1+(1+\varphi)c^r+(1+\varphi)(1+\varepsilon_2)\sum_{i=1}^{r-1}c^i\right)a^2x^2 \\ &\quad+(1+\varepsilon_1^{-1})(x+y_1)^2+c^{-1}(1+\varphi)(1+\varepsilon_2^{-1})\sum_{i=2}^{r}c^i\left(x+y_1+\ldots+y_i\right)^2\end{aligned} \tag{2.28}$$

It follows from (2.28) and (2.22), that there exists $\sigma \in [0,1)$ such that $\bar{V}(x+dx+y_1,y_2,...,y_r,u) \le \sigma \bar{V}(x,y_1,y_2,...,y_r)$ holds for all $(x,y_1,...,y_r) \in \Re^{r+1}$ and $|d| \le a$, provided that there exist $\varepsilon_1 > 0$ and $\varepsilon_2 > 0$ so that the following inequalities hold:

$$\begin{aligned} &\left(1+\varepsilon_1+(1+\varphi)c^r+(1+\varphi)(1+\varepsilon_2)\sum_{i=1}^{r-1} c^i\right)a^2 < 1 \\ &1+\varepsilon_1^{-1} < c(1+\varphi) \\ &1+\varepsilon_2^{-1} < c \end{aligned} \tag{2.29}$$

Inequalities (2.29) hold for appropriate $\varepsilon_1 > 0$ and $\varepsilon_2 > 0$ provided that the following inequality holds:

$$a^2 < \frac{s}{1+s\left(1+\dfrac{c^{r+1}-c^r+c^{r-1}-c}{(c-1)^2}\right)+s^2\left(\dfrac{c^{r+1}-c^r+c^{r-1}-c}{(c-1)^2}\right)} \tag{2.30}$$

where $s = c(1+\varphi)-1 > 0$. The greatest value for the fraction $\frac{s}{1+s\left(1+\frac{c^{r+1}-c^r+c^{r-1}-c}{(c-1)^2}\right)+s^2\left(\frac{c^{r+1}-c^r+c^{r-1}-c}{(c-1)^2}\right)}$ is obtained for $s = \frac{c-1}{\sqrt{c^{r+1}-c^r+c^{r-1}-c}}$. The value of $c > 1$, which maximizes the right hand side of (2.30) for $s = \frac{c-1}{\sqrt{c^{r+1}-c^r+c^{r-1}-c}}$ can be found numerically.

What have we found so far? We have shown that $0 \in \Re \times \Re^r$ is robustly globally asymptotically stable for the closed-loop system (2.18) with (2.20) provided that $a < A_r$. The value of $a < A_r$ is estimated by the necessary condition (2.21) and the sufficient condition (2.30) with $s = \frac{c-1}{\sqrt{c^{r+1}-c^r+c^{r-1}-c}}$. The results are shown in Table 1.

| $r$ | $A_r$ |
|---|---|
| 0 | $A_0 = 1$ |
| 1 | $A_1 = 0.5$ |
| 2 | $A_2 \in [0.3311, 0.3333]$ |
| 3 | $A_3 \in [0.2451, 0.25]$ |
| 4 | $A_4 \in [0.1923, 0.2]$ |
| 5 | $A_5 \in [0.1573, 0.1667]$ |
| 6 | $A_6 \in [0.1326, 0.1429]$ |
| 7 | $A_7 \in [0.1144, 0.125]$ |
| 8 | $A_8 \in [0.1005, 0.1112]$ |
| 9 | $A_9 \in [0.0896, 0.1]$ |
| 10 | $A_{10} \in [0.0807, 0.0909]$ |
| 15 | $A_{15} \in [0.0539, 0.0625]$ |
| 20 | $A_{20} \in [0.0404, 0.0476]$ |

**Table 1:** Results of the robustness analysis for system (2.18) with (2.20)

It is clear that as the value of the input delay $r$ increases, the value of $A_r$ with the property that $0 \in \Re \times \Re^r$ is robustly globally asymptotically stable for the closed-loop system (2.18) with (2.20), $|d| \leq a$ and $a < A_r$, decreases rapidly. In other words, as the value of the input delay $r$ increases, the sensitivity with respect to uncertain parameter $d \in \Re$ is magnified.

The example will be studied further. $\triangleleft$

The results of Example 2.3 are expected. However, it should be emphasized that the results are not discouraging for the use of predictor feedback: the sensitivity with respect to modeling errors is not magnified because of the use of the predictor feedback. This happens because the control problem itself is difficult when the value of the input delay increases. On the other hand, we should seek predictor-based feedback laws that minimize the sensitivity with respect to modeling errors as much as possible. This is the topic of the following section.

## 3. Lyapunov Redesign

As remarked in the previous section, it is very important to design a feedback law that minimizes the sensitivity with respect to modeling errors. It is important to emphasize that although the "nominal feedback law" $u = k(x)$ may be "optimal" in the sense that minimizes some measure of the sensitivity of the corresponding closed-loop system with respect to modeling errors, this is not necessarily true for the feedback law proposed by Lemma 2.1 or Lemma 2.2.

In order to design a feedback law that minimizes the sensitivity with respect to modeling errors, we exploit the Lyapunov function proposed by Lemma 2.1 and Lemma 2.2. We consider the problem of robust global feedback stabilization of the equilibrium point $0 \in \Re^n \times \Re^{mr}$ for the uncertain control system (1.5).

The procedure that we are proposing includes the following steps:

The Lyapunov redesign procedure:

Step 1: Find a value $d_0 \in D$ for the uncertainty $d \in D$, which is nominal in a certain sense. Define the "nominal" vector field $f(x,u) = F(d_0, x, u)$, for all $(x,u) \in \Re^n \times \Re^m$.

Step 2: Verify assumption (H1) for the corresponding system (1.1) with $r = 0$. More specifically, find a continuous function $k : \Re^n \to \Re^m$ with $k(0) = 0$, a continuous, positive definite and radially unbounded function $V : \Re^n \to \Re_+$ and a constant $\lambda \in [0,1)$ such that inequality (2.1) holds.

Step 3: Use a family of functions $a_i \in K_\infty$ ($i = 1,...,r$) with $a_i(s) \leq a_{i+1}(s)$ for all $i = 1,...,r-1$ and $s \geq 0$, and a constant $c > 1$ and define the function $\bar{V} : \Re^n \times \Re^{rm} \to \Re_+$ by means of (2.11).

Step 4: For each $(x, y_1, ..., y_r) \in \Re^n \times \Re^{rm}$ solve the minimax problem

$$\min_{y_{r+1} \in \Re^m} \max_{d \in D} \bar{V}(F(d, x, y_1), y_2, ..., y_{r+1}) \tag{3.1}$$

If problem (3.1) is solvable for every $(x, y_1, ..., y_r) \in \Re^n \times \Re^{rm}$ and if the following inequality holds for all $(x, y_1, ..., y_r) \in \Re^n \times \Re^{rm}$ :

$$\min_{y_{r+1} \in \Re^m} \max_{d \in D} \bar{V}(F(d, x, y_1), y_2, ..., y_{r+1}) \leq \bar{V}(x, y_1, ..., y_r) - \rho\left(\bar{V}(x, y_1, ..., y_r)\right) \tag{3.2}$$

for certain continuous and positive definite function $\rho : \Re_+ \to \Re_+$ then the robust global feedback stabilizer $K(x, y_1, ..., y_r)$ can be defined as any of the minimizers of the minimax problem (3.1), i.e., the robust global feedback stabilizer $K(x, y_1, ..., y_r)$ satisfies for all $(x, y_1, ..., y_r) \in \Re^n \times \Re^{rm}$

$$\min_{y_{r+1} \in \Re^m} \max_{d \in D} \bar{V}(F(d, x, y_1), y_2, ..., y_{r+1}) = \max_{d \in D} \bar{V}(F(d, x, y_1), y_2, ..., y_r, K(x, y_1, ..., y_r)) \tag{3.3}$$

The procedure that we just described has many "open issues":

1) Under what conditions will the minimax problem (3.1) be solvable for all $(x, y_1, ..., y_r) \in \Re^n \times \Re^{rm}$ ?
2) Under what conditions does there exist a continuous and positive definite function $\rho : \Re_+ \to \Re_+$ such that (3.2) holds?
3) What are the regularity properties for the function $K(x, y_1, ..., y_r)$ that satisfies (3.3)?
4) How can we select the nominal value for the uncertainty parameter $d_0 \in D$, the constant $c > 1$ and the family of functions $a_i \in K_\infty$ $(i = 1, ..., r)$ with $a_i(s) \leq a_{i+1}(s)$ for all $i = 1, ..., r-1$ and $s \geq 0$?

The following theorem answers questions (1) and (3) above.

**Theorem 3.1:** *Assume that* $D \subseteq \Re^l$ *is a non-empty, compact set. Then the minimax problem (3.1) is solvable for all* $(x, y_1, ..., y_r) \in \Re^n \times \Re^{rm}$ *and every function* $K : \Re^n \times \Re^{rm} \to \Re^m$, *that satisfies (3.3) for all* $(x, y_1, ..., y_r) \in \Re^n \times \Re^{rm}$, *is locally bounded. Moreover, there exists a measurable function* $K : \Re^n \times \Re^{rm} \to \Re^m$ *that satisfies (3.3) for all* $(x, y_1, ..., y_r) \in \Re^n \times \Re^{rm}$. *Finally, if there exists an open set* $O \subseteq \Re^n \times \Re^{rm}$, *such that the minimax problem (3.1) has a unique solution for all* $(x, y_1, ..., y_r) \in O$, *then every function* $K : \Re^n \times \Re^{rm} \to \Re^m$, *that satisfies (3.3) for all* $(x, y_1, ..., y_r) \in \Re^n \times \Re^{rm}$, *is continuous on* $O \subseteq \Re^n \times \Re^{rm}$.

**Proof:** Since $\bar{V} : \Re^n \times \Re^{rm} \to \Re_+$ is a continuous, positive definite and radially unbounded function, by virtue of Lemma 3.5 in [15], there exist functions $a_1, a_2 \in K_\infty$ such that:

$$a_1\left(|x, y_1, y_2, ..., y_r|\right) \leq \bar{V}(x, y_1, y_2, ..., y_r) \leq a_2\left(|x, y_1, y_2, ..., y_r|\right), \text{ for all } (x, y_1, ..., y_r) \in \Re^n \times \Re^{rm} \tag{3.4}$$

Define for all $(x, y_1, ..., y_r, y_{r+1}) \in \Re^n \times \Re^{(r+1)m}$ :

$$\Psi(x, y_1, y_2, ..., y_{r+1}) := \max_{d \in D} \bar{V}(F(d, x, y_1), y_2, ..., y_{r+1}) \tag{3.5}$$

Theorem 1.4.16 in [2], in conjunction with continuity of $\bar{V} : \Re^n \times \Re^{rm} \to \Re_+$, $F : D \times \Re^n \times \Re^m \to \Re^n$ and compactness of $D \subseteq \Re^l$ implies that $\Psi : \Re^n \times \Re^{(r+1)m} \to \Re_+$ as defined by (3.5) is continuous.

Since the mapping $\Psi:\Re^n\times\Re^{(r+1)m}\to\Re_+$ is bounded from below, we can define for all $(x,y_1,...,y_r)\in\Re^n\times\Re^{rm}$:

$$\tilde{V}(x,y_1,...,y_r):=\inf\left\{\Psi(x,y_1,...,y_r,y_{r+1})\,;\,y_{r+1}\in\Re^m\right\} \tag{3.6}$$

$$M(x,y_1,...,y_r):=\left\{y_{r+1}\in\Re^m:\tilde{V}(x,y_1,...,y_r)=\Psi(x,,y_1,...,y_r,y_{r+1})\right\} \tag{3.7}$$

In order to show that every function $K:\Re^n\times\Re^{rm}\to\Re^m$, that satisfies (3.3) for all $(x,y_1,...,y_r)\in\Re^n\times\Re^{rm}$, is locally bounded, it suffices to show that the set-valued mapping $M(x,y_1,...,y_r)\subseteq\Re^m$ defined by (3.7) is non-empty for all $(x,y_1,...,y_r)\in\Re^n\times\Re^{rm}$ and locally bounded. Define:

$$p(x,y_1,...,y_r):=a_1^{-1}\left(\Psi(x,y_1,...,y_r,0)+1\right) \tag{3.8}$$

and notice that the mapping $p(x,y_1,...,y_r)$ is a continuous, positive function. Definitions (3.5), (3.6), (3.8) and the left hand side inequality (3.4) imply that for each fixed $(x,y_1,...,y_r)\in\Re^n\times\Re^{rm}$ we have:

$$\begin{aligned}
&\tilde{V}(x,y_1,...,y_r)=\\
&\min\left(\inf\left\{\Psi(x,y_1,...,y_r,y_{r+1})\,;\,|y_{r+1}|\le p(x,y_1,...,y_r)\right\},\inf\left\{\Psi(x,y_1,...,y_r,y_{r+1})\,;\,|y_{r+1}|>p(x,y_1,...,y_r)\right\}\right)\\
&\ge\min\left(\inf\left\{\Psi(x,y_1,...,y_r,y_{r+1})\,;\,|y_{r+1}|\le p(x,y_1,...,y_r)\right\},\inf\left\{a_1(|y_{r+1}|)\,;\,|y_{r+1}|>p(x,y_1,...,y_r)\right\}\right)\\
&\ge\min\left(\inf\left\{\Psi(x,y_1,...,y_r,y_{r+1})\,;\,|y_{r+1}|\le p(x,y_1,...,y_r)\right\},\ \Psi(x,y_1,...,y_r,0)+1\right)
\end{aligned}$$

Clearly, since $\tilde{V}(x,y_1,...,y_r)\le\Psi(x,y_1,...,y_r,y_{r+1})$, the above inequality implies that the case $\min\left(\inf\left\{\Psi(x,y_1,...,y_r,y_{r+1})\,;\,|y_{r+1}|\le p(x,y_1,...,y_r)\right\},\ \Psi(x,y_1,...,y_r,0)+1\right)=\Psi(x,y_1,...,y_r,0)+1$ cannot happen. Thus we conclude that:

$$\tilde{V}(x,y_1,...,y_r)=\inf\left\{\Psi(x,y_1,...,y_r,y_{r+1})\,;\,|y_{r+1}|\le p(x,y_1,...,y_r)\right\} \tag{3.9}$$

Equality (3.9) in conjunction with continuity of $\Psi(x,y_1,y_2,...,y_{r+1})$ implies that the set-valued map $M(x,y_1,...,y_r)\subseteq\Re^m$, as defined by (3.7), is non-empty for each $(x,y_1,...,y_r)\in\Re^n\times\Re^{rm}$. Continuity of the mapping $p(x,y_1,...,y_r)$ and definitions (3.6), (3.7) in conjunction with (3.9) imply that the set-valued map $M(x,y_1,...,y_r)\subseteq\Re^m$ is locally bounded: notice that every $u\in M(x,y_1,...,y_r)$ satisfies $|u|\le p(x,y_1,...,y_r)$.

Moreover, continuity of the mapping $p(x,y_1,...,y_r)$, Corollary 1.4.10 in [2] (and the remark just after the statement of Corollary 1.4.20 in [2], page 43), Theorem 1.4.16 in [2] (page 48) and equality (3.9) imply that the mapping $(x,y_1,...,y_r)\to\tilde{V}(x,y_1,...,y_r)$ is continuous. Continuity of the mappings $(x,y_1,...,y_r)\to\tilde{V}(x,y_1,...,y_r)$ and $\Psi:\Re^n\times\Re^{(r+1)m}\to\Re_+$ in conjunction with definition (3.7) and statement (c) on page 150 in [3] imply that the set-valued mapping $M(x,y_1,...,y_r)$ is measurable. Consequently, Theorem 5.3 on page 151 in [3] implies that there exists a measurable function $K:\Re^n\times\Re^{rm}\to\Re^m$ that satisfies $K(x,y_1,...,y_r)\in M(x,y_1,...,y_r)$ for all $(x,y_1,...,y_r)\in\Re^n\times\Re^{rm}$. Therefore, definitions (3.5), (3.6) and (3.7) imply that $K:\Re^n\times\Re^{rm}\to\Re^m$ satisfies (3.3) for all $(x,y_1,...,y_r)\in\Re^n\times\Re^{rm}$.

In order to show that the last assertion of the theorem, it suffices to show that the set-valued map $\mathcal{M}(x, y_1, ..., y_r) \subseteq \Re^m$ is upper semi-continuous. Indeed, this automatically implies that if $\mathcal{M}(x, y_1, ..., y_r) \subseteq \Re^m$ is a singleton for all $(x, y_1, ..., y_r) \in O$, where $O \subseteq \Re^n \times \Re^{rm}$ is an open set, i.e., $\mathcal{M}(x, y_1, ..., y_r) = \{\varphi(x, y_1, ..., y_r)\}$, then $\varphi(x, y_1, ..., y_r)$ is continuous on $O \subseteq \Re^n \times \Re^{rm}$.

In order to show that $\mathcal{M}(x, y_1, ..., y_r) \subseteq \Re^m$ is upper semi-continuous, it suffices to prove that for every $z_r = (x, y_1, ..., y_r) \in \Re^n \times \Re^{rm}$ and for every $\varepsilon > 0$ there exists $\delta > 0$ such that

$$|w - z_r| < \delta \quad \Rightarrow \mathcal{M}(w) \subset \mathcal{M}(z_r) + \varepsilon B$$

The proof is made by contradiction. Suppose the contrary: there exists $z_r = (x, y_1, ..., y_r) \in \Re^n \times \Re^{rm}$ and $\varepsilon > 0$, such that for all $\delta > 0$, there exists $w \in \{z_r\} + \delta B$ and $u' \in \mathcal{M}(w)$ with $|u' - u| \geq \varepsilon$, for all $u \in \mathcal{M}(z_r)$. Clearly, this implies the existence of a sequence $\{(w_j, u'_j)\}_{j=1}^{\infty}$ with $w_j \to z_r$, $u'_j \in \mathcal{M}(w_j)$ and $|u'_j - u| \geq \varepsilon$, for all $u \in \mathcal{M}(z_r)$ and $j = 1,2,...$. On the other hand since $u'_j$ is bounded, it contains a convergent subsequence $u'_i \to \bar{u} \notin \mathcal{M}(z_r)$. By continuity of the mappings $\tilde{V}(x, y_1, ..., y_r)$ and $\Psi(x, y_1, ..., y_r, u)$, we have: $\tilde{V}(w_i) \to \tilde{V}(z_r)$ and $\tilde{V}(w_i) = \Psi(w_i, u'_i) \to \Psi(z_r, \bar{u})$. Consequently, we must have: $\tilde{V}(z_r) = \Psi(z_r, \bar{u})$, which, by virtue of definition (3.7) implies that $\bar{u} \in \mathcal{M}(z_r)$, a contradiction.

The proof is complete. $\triangleleft$

The answer to question (4) above is an open problem which is directly related to the answer to question (2). However, there are some simple cases for which we can give explicit formulae for the stabilizing feedback. Next we consider the simple single-input case (1.2) with $F(d, x, u) = Ax + Bu + dGx$, $x \in \Re^n, u \in \Re, d \in D = [-a, a] \subseteq \Re$. The reader should notice that even in this "almost linear" case, the proposed feedback is nonlinear: it is a homogeneous function of degree 1.

**Theorem 3.2:** *Consider the single input discrete-time system*

$$\begin{aligned} &x(t+1) = Ax(t) + Bu(t) + d(t)Gx(t) \\ &x(t) \in \Re^n, u(t) \in \Re, d(t) \in D = [-a, a] \subseteq \Re \end{aligned} \tag{3.10}$$

*where* $A \in \Re^{n \times n}, G \in \Re^{n \times n}$ *are constant matrices,* $a \geq 0$ *is a constant and* $B \in \Re^n$ *is a constant vector. Suppose that there is a vector* $k \in \Re^n$*, a constant* $\lambda \in [0,1)$ *and a symmetric positive definite matrix* $P \in \Re^{n \times n}$ *such that the following inequality holds for all* $x \in \Re^n$*:*

$$x'(A + Bk')'P(A + Bk')x \leq \lambda x'Px \tag{3.11}$$

*Let* $r \geq 2$ *be a positive integer, and let* $c > 1$*,* $\varphi > 0$*,* $\sigma \in [0,1)$ *be constants. Define:*

$$p := c^r (B'PB + \varphi) \tag{3.12}$$

$$L(x) := c^r (B'PA - \varphi k') A^{r-1} Gx \tag{3.13}$$

$$\kappa(x,y_1,...,y_r) := \left(Ax+By_1\right)' PGx + \sum_{i=1}^{r-1} c^i y_{i+1}\left(B'PA-\varphi k'\right)A^{i-1}Gx$$
$$+\sum_{i=1}^{r} c^i \left( A^i x + \sum_{j=1}^{i} A^{i-j} B y_j \right)' \left(A'PA+\varphi kk'\right)A^{i-1}Gx \tag{3.14}$$

$$b(x,y_1,...,y_r) := c^r \left(B'PA-\varphi k'\right)\left( A^r x + \sum_{j=1}^{r} A^{r-j} B y_j \right) \tag{3.15}$$

$$c(x,y_1,...,y_r) := a^2 \sum_{i=0}^{r} c^i x'G'(A^i)'P(A^i)Gx + a^2\varphi \sum_{i=1}^{r} c^i \left(k'A^{i-1}Gx\right)^2$$
$$+(1-\sigma c)\sum_{i=1}^{r} c^{i-1}\left( A^i x + \sum_{j=1}^{i} A^{i-j} B y_j \right)' P \left( A^i x + \sum_{j=1}^{i} A^{i-j} B y_j \right)$$
$$+c^r \left( A^r x + \sum_{j=1}^{r} A^{r-j} B y_j \right)' \left(A'PA+\varphi kk'\right)\left( A^r x + \sum_{j=1}^{r} A^{r-j} B y_j \right) \tag{3.16}$$
$$+(1-\sigma c)\varphi \sum_{i=2}^{r} c^{i-1}\left( y_i - k'A^{i-1}x - k'\sum_{j=1}^{i-1} A^{i-1-j} B y_j \right)^2 - \sigma\, x'Px - \sigma c \varphi \left(y_1 - k'x\right)^2$$

*Consider the continuous, homogeneous of degree 1, function defined by:*

$$K(x,y_1,...,y_r) := \begin{cases} -L^{-1}(x)\kappa(x,y_1,...,y_r) & \text{if} \quad \left|p\kappa(x,y_1,...,y_r) - b(x,y_1,...,y_r)L(x)\right| < aL^2(x) \\ -p^{-1}\left(aL(x)+b(x,y_1,...,y_r)\right) & \text{if} \quad p\kappa(x,y_1,...,y_r) - b(x,y_1,...,y_r)L(x) \geq aL^2(x) \\ p^{-1}\left(aL(x)-b(x,y_1,...,y_r)\right) & \text{if} \quad p\kappa(x,y_1,...,y_r) - b(x,y_1,...,y_r)L(x) \leq -aL^2(x) \end{cases} \tag{3.17}$$

*Suppose that the following inequalities hold:*

$$p\left(\frac{\kappa(x,y_1,...,y_r)}{L(x)}\right)^2 - 2b(x,y_1,...,y_r)\frac{\kappa(x,y_1,...,y_r)}{L(x)} + c(x,y_1,...,y_r) \leq 0,$$
$$\text{for all } (x,y_1,...,y_r) \in \Re^n \times \Re^r \text{ with } \left|p\kappa(x,y_1,...,y_r) - b(x,y_1,...,y_r)L(x)\right| < aL^2(x) \tag{3.18}$$

$$-p^{-1}\left(aL(x)+b(x,y_1,...,y_r)\right)^2 + c(x,y_1,...,y_r) + 2a\kappa(x,y_1,...,y_r) \leq 0,$$
$$\text{for all } (x,y_1,...,y_r) \in \Re^n \times \Re^r \text{ with } p\kappa(x,y_1,...,y_r) - b(x,y_1,...,y_r)L(x) \geq aL^2(x) \tag{3.19}$$

$$-p^{-1}\left(aL(x)-b(x,y_1,...,y_r)\right)^2 + c(x,y_1,...,y_r) - 2a\kappa(x,y_1,...,y_r) \leq 0,$$
$$\text{for all } (x,y_1,...,y_r) \in \Re^n \times \Re^r \text{ with } p\kappa(x,y_1,...,y_r) - b(x,y_1,...,y_r)L(x) \leq -aL^2(x) \tag{3.20}$$

*Then* $0 \in \Re^n \times \Re^r$ *is robustly globally exponentially stable for the closed-loop system:*

$$\begin{aligned} & x(t+1) = Ax(t) + By_1(t) + d(t)Gx(t) \\ & y_i(t+1) = y_{i+1}(t) \quad , \quad i = 1,...,r \\ & x(t) \in \Re^n \, , y_i(t) \in \Re \left(i = 1,...,r\right), \\ & y_{r+1}(t) = u(t) \in \Re \, , d(t) \in D = [-a,a] \end{aligned} \tag{3.21}$$

*with (1.7) and (3.17).*

**Remark 3.3:** Inequalities (3.18), (3.19) and (3.20) cannot be guaranteed by design. They are the inequalities which allow us to calculate (numerically) the lowest upper bound for the uncertainty magnitude $a>0$ for which robust global exponential stability holds for the closed-loop system (3.21) with (1.7) and (3.17). The use of inequalities (3.18), (3.19) and (3.20) is illustrated in the example below.

**Proof:** For $c>1$, $\varphi>0$, we consider the Lyapunov function $\overline{V}:\Re^n \times \Re^r \to \Re_+$ defined by (2.11) with $a_i(s) \equiv \varphi s^2$ for $i=1,...,r$, $k(x)=k'x$ and $V(x)=x'Px$ for the disturbance-free discrete-time system (1.4) with $f(x,u)=Ax+Bu$, $x\in\Re^n, u\in\Re$ (this is the control system that corresponds to the nominal value of the disturbance $d=0$):

$$\begin{aligned}\overline{V}(z_r) &:= x'Px+\sum_{i=1}^{r} c^i \left( A^i x+\sum_{j=1}^{i} A^{i-j}By_j \right)' P \left( A^i x+\sum_{j=1}^{i} A^{i-j}By_j \right) \\ &+c\varphi(y_1-k'x)^2+\varphi\sum_{i=2}^{r} c^i \left( y_i-k'A^{i-1}x-k'\sum_{j=1}^{i-1} A^{i-1-j}By_j \right)^2\end{aligned} \tag{3.22}$$

where $z_r:=(x,y_1,...,y_r)\in\Re^n\times\Re^r$. Notice that the formula (3.22) coincides with formula (2.11), since we have $F_i(z_i)=A^i x+\sum_{j=1}^{i} A^{i-j}By_j$, where the vector fields $F_i(z_i)$, $i=0,1,...,r+1$ are defined by the recursive formula (2.9). Using (3.22) and definitions (3.12), (3.13), (3.14), (3.15), (3.16), we obtain for all $(x,y_1,...,y_r,u,d)\in\Re^n\times\Re^r\times\Re\times\Re$:

$$\begin{aligned}&\overline{V}(Ax+By_1+dGx,y_2,...,y_r,u)= \\ &pu^2+2b(x,y_1,y_2,...,y_r)u+2d\big(\kappa(x,y_1,y_2,...,y_r)+L(x)u\big) \\ &+c(x,y_1,y_2,...,y_r)+\sigma\overline{V}(x,y_1,y_2,...,y_r) \\ &+\left(d^2-a^2\right)\left(\sum_{i=0}^{r} c^i x'G'(A^i)'P(A^i)Gx+\varphi\sum_{i=1}^{r} c^i\left(dk'A^{i-1}Gx\right)^2\right)\end{aligned} \tag{3.23}$$

Since $P\in\Re^{n\times n}$ is positive definite, it follows that the coefficient of $d^2$ in (3.23) is non-negative, i.e. $\sum_{i=0}^{r} c^i x'G'(A^i)'P(A^i)Gx+\varphi\sum_{i=1}^{r} c^i\left(dk'A^{i-1}Gx\right)^2\geq 0$. Therefore, it follows from (3.23) for all $(x,y_1,...,y_r,u)\in\Re^n\times\Re^r\times\Re$:

$$\begin{aligned}&\max_{|d|\leq a}\overline{V}(Ax+By_1+dGx,y_2,...,y_r,u)= \\ &pu^2+2b(x,y_1,y_2,...,y_r)u+2a\left|\kappa(x,y_1,y_2,...,y_r)+L(x)u\right| \\ &+c(x,y_1,y_2,...,y_r)+\sigma\overline{V}(x,y_1,y_2,...,y_r)\end{aligned} \tag{3.24}$$

It follows from the minimization of the function defined in (3.24) that the minimizer must satisfy $u=K(x,y_1,...,y_r)$, where $K:\Re^n\times\Re^r\to\Re$ is defined by (3.17). Finally, inequalities (3.18), (3.19) and (3.20) guarantee that the inequality

$$\max_{|d|\leq a}\overline{V}(Ax+By_1+dGx,y_2,...,y_r,K(x,y_1,y_2,...,y_r))\leq\sigma\overline{V}(x,y_1,y_2,...,y_r) \tag{3.25}$$

holds for all $(x,y_1,...,y_r)\in\Re^n\times\Re^r$. The conclusion of the theorem is a consequence of (3.25) and Proposition 2.3 in [12]. The proof is complete. ◁

A similar result with that of Theorem 3.2 holds for the case $r = 1$.

Notice that Theorem 3.2 does not give a complete answer to the Lyapunov redesign procedure for the simple system (3.21). However, the control practitioner can use the formulae in the statement of Theorem 3.2 and select values for the constants $c > 1$, $\varphi > 0$, $\sigma \in [0,1)$ so that the value of $a \geq 0$ becomes as large as possible and thus minimize the sensitivity with respect to modeling errors.

**Example 3.4:** The importance of Lyapunov redesign will be illustrated by means of system (2.18) for $r = 1$, which was studied in Example 2.3. We consider again we consider the Lyapunov function defined by (2.22) for $c > 0$ and $\varphi > -1$. For $r = 1$, we get for all $(x, y_1) \in \Re^2$, $d \in \Re$:

$$\begin{aligned} \bar{V}(x + dx + y_1, u) &= (x + dx + y_1)^2 + (1+\varphi)c(x + dx + y_1 + u)^2 \\ &= (x + y_1)^2 + (1+\varphi)c(x + y_1 + u)^2 + (1 + c(1+\varphi))d^2 x^2 + 2dx(2x + 2y_1 + u) \end{aligned}$$

The above equality implies that:

$$\begin{aligned} \max_{|d| \leq a} \bar{V}(x + dx + y_1, u) = {} & qu^2 + (1+q)a^2 x^2 + 2a\left|2x^2 + 2xy_1 + xu\right| \\ & + (1+q)(x + y_1)^2 + 2qu(x + y_1) \end{aligned}$$

where $q = c(1+\varphi) > 0$. The feedback law is defined as the minimizer of the above quantity, i.e.,

$$u = -\left(1 + \frac{a}{q}\right)x - y_1, \text{ for } x^2 + xy_1 \geq \frac{a}{q}x^2 \tag{3.26}$$

$$u = -\left(1 - \frac{a}{q}\right)x - y_1, \text{ for } x^2 + xy_1 \leq -\frac{a}{q}x^2 \tag{3.27}$$

$$u = -2x - 2y_1, \text{ for } |x + y_1| < \frac{a}{q}|x| \tag{3.28}$$

Figure 1 shows the three regions in the state space which are involved in (3.26), (3.27) and (3.28). Notice that the feedback law defined by (3.26), (3.27) and (3.28) is a continuous, piecewise linear feedback law, which is homogeneous of degree 1.

The reader should notice the difference between the above feedback law and the feedback law defined by (2.20) with $r = 1$, which was obtained with no Lyapunov redesign. In order, to find the value for $q = c(1+\varphi) > 0$ that allows $a \geq 0$ to be as large as possible, we follow a robustness analysis similar to the analysis of Example 2.3. The existence of $\sigma \in [0,1)$ so that $\max_{|d| \leq a} \bar{V}(x + dx + y_1, u) \leq \sigma \bar{V}(x, y_1)$ for all $(x, y_1) \in \Re^2$ and $u$ given by (3.26), (3.27) and (3.38) is equivalent to the following inequalities:

$$\left(2a-\frac{a^2}{q}+(1+q)a^2-\sigma+1-\sigma q\right)x^2+2(a+1-\sigma q)xy_1+(1-\sigma q)y_1^2\le 0\text{, for } x^2+xy_1\ge\frac{a}{q}x^2$$

$$\left((1+q)a^2-2a-\frac{a^2}{q}-\sigma+1-\sigma q\right)x^2-2(a+\sigma q-1)xy_1+(1-\sigma q)y_1^2\le 0\text{, for } x^2+xy_1\le-\frac{a}{q}x^2$$

$$(1+q-\sigma q)(x+y_1)^2+\left((1+q)a^2-\sigma\right)x^2\le 0\text{, for } |x+y_1|<\frac{a}{q}|x|$$

The existence of $\sigma\in[0,1)$ that satisfies the above inequalities is equivalent to the following inequalities:

$$\left(2a-\frac{a^2}{q}+(1+q)a^2-1\right)\cos^2(\theta)+(a+1-q)\sin(2\theta)<q-1\text{, for all } \theta\in[0,2\pi) \text{ with } \sin(2\theta)\ge 2\left(\frac{a}{q}-1\right)\cos^2(\theta)$$

$$\left((1+q)a^2-2a-\frac{a^2}{q}-1\right)\cos^2(\theta)+(1-a-q)\sin(2\theta)<q-1\text{, for all } \theta\in[0,2\pi) \text{ with } \sin(2\theta)\le -2\left(\frac{a}{q}+1\right)\cos^2(\theta)$$

$$\left((1+q)a^2-1\right)\cos^2(\theta)+1+\sin(2\theta)<0\text{, for all } \theta\in[0,2\pi) \text{ with } -2\left(\frac{a}{q}+1\right)\cos^2(\theta)<\sin(2\theta)<2\left(\frac{a}{q}-1\right)\cos^2(\theta)$$

The numerical evaluation of all the above quantities shows that for $a=0.535$ all the above inequalities hold with $q=1.81$. The reader should notice the improvement compared to the feedback design with no Lyapunov redesign of Example 2.3, where the necessary and sufficient condition for robust global asymptotic stability was $a<0.5$.

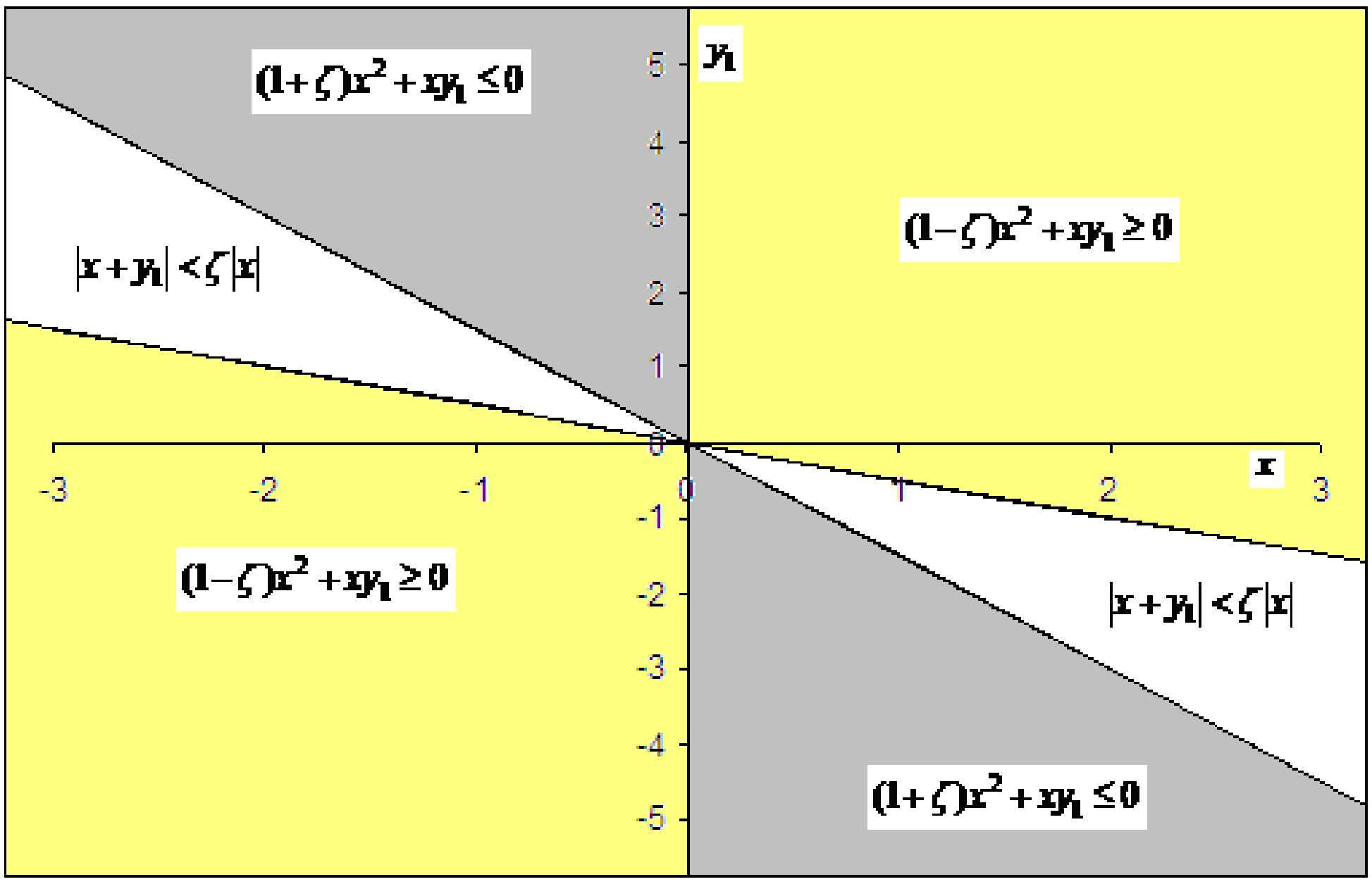

**Figure 1:** The three regions in the state space, which are involved in (3.26), (3.27) and (3.28), for $\zeta=\frac{a}{q}=0.5$.

Even better results can be obtained if we notice that a feedback stabilizer for the delay-free system can be given by the formula $k(x)=-\beta x$, where $\beta\in(0,2)$. ◁

## 4. Concluding Remarks

This work studied the design problem of feedback stabilizers for discrete-time systems with input delays. A backstepping procedure is proposed for disturbance-free discrete-time systems. The feedback law designed by using backstepping coincides with the predictor-based feedback law used in continuous-time systems with input delays. As in the continuous-time case, the backstepping procedure allows a simultaneous determination of a stabilizing feedback and a Lyapunov function for the corresponding closed-loop system (Lemma 2.1 and Lemma 2.2).

However, simple examples demonstrate that the sensitivity of the closed-loop system with respect to modeling errors increases as the value of the delay increases (Example 2.3). The paper proposed a Lyapunov redesign procedure which can minimize the effect of the uncertainty. The point values of the feedback law can be found as the solution of a specific minimax problem (Theorem 3.1) Specific results are provided for linear single-input discrete-time systems with multiplicative uncertainty (Theorem 3.2). The feedback law that guarantees robust global exponential stability is a nonlinear, homogeneous of degree 1 feedback law. Explicit inequalities allow the determination of the lowest upper bound for the magnitude of the uncertainty for which robust global exponential stability holds for the closed-loop system. A simple example showed the importance of the Lyapunov redesign procedure (Example 3.4).

## References

[1] Artstein, Z., "Linear systems with delayed controls: A reduction", *IEEE Transactions on Automatic Control*, 27, 1982, 869-879.

[2] Aubin, J. P. and H. Frankowska, *Set-Valued Analysis*, Birkhauser, Boston, 1990.

[3] Clarke, F. H., Yu. S. Ledyaev, R.J. Stern and P.R. Wolenski, *Nonsmooth Analysis and Control Theory*, Spinger-Verlag, New York, 1998.

[4] Garcia, P., A. Gonzalez, P. Castillo, R. Lozano and P. Albertos, "Robustness of a discrete-time predictor-based controller for time-varying measurement delay", *Proceedings of the 9th IFAC workshop on time delay systems*, Prague, Czech Republic, 2010.

[5] Gonzalez, A., A. Sala, P. Albertos, "Predictor-based stabilization of discrete time-varying input-delay systems", *Automatica* 48, 2012, 454–457.

[6] Gonzalez, A., A. Sala, P. Garcia and P. Albertos, "Robustness analysis of discrete predictor-based controllers for input-delay systems", *International Journal of Systems Science*, 44(2), 2013, 232-239.

[7] Grüne, L., and K. Worthmann "Nonlinear sampled-data redesign: analytical formulas and their practical implementation", *Proceedings of the 18th International Symposium on Mathematical Theory of Networks and Systems (MTNS2008),* CD-Rom, Paper 084.pdf, Blacksburg, Virginia, 2008.

[8] Grüne, L., K. Worthmann and D. Nešic, "Continuous-time controller redesign for digital implementation: a trajectory based approach", *Automatica,* 44, 2008, 225 - 232.

[9] Henson, M. A. and D. E. Seborg, "Time delay compensation for nonlinear processes", *Industrial and Engineering Chemistry Research*, 33(6), 1994, 1493–1500.

[10] Jiang, Z.P. and Y. Wang, "A Converse Lyapunov Theorem for Discrete-Time Systems with Disturbances", *Systems and Control Letters*, 45(1), 2002, 49-58.

[11] Karafyllis, I., and S. Kotsios, “Necessary and Sufficient Conditions for Robust Global Asymptotic Stabilization of Discrete-Time Systems”, *Journal of Difference Equations and Applications*, 12(7), 2006, pp. 741-768.

[12] Karafyllis, I., and Z.-P. Jiang, *Stability and Stabilization of Nonlinear Systems*, Springer-Verlag London (Series: Communications and Control Engineering), 2011.

[13] Karafyllis, I. and M. Krstic, “Nonlinear Stabilization under Sampled and Delayed Measurements, and with Inputs Subject to Delay and Zero-Order Hold”, *IEEE Transactions on Automatic Control*, 57(5), 2012, 1141-1154.

[14] Kellett, C. M. and A.R. Teel, “Smooth Lyapunov Functions and Robustness of Stability for Difference Inclusions”, *Systems and Control Letters*, 52(5), 2004, 395-405.

[15] Khalil, H. K., *Nonlinear Systems*, 2nd Edition, Prentice-Hall, 1996.

[16] Krstic, M., I. Kanellakopoulos, and P. V. Kokotovic, *Nonlinear and Adaptive Control Design*, Wiley, 1995.

[17] Krstic, M., *Delay Compensation for Nonlinear, Adaptive, and PDE Systems*, Birkhäuser, Boston, 2009.

[18] Krstic, M., “Input Delay Compensation for Forward Complete and Strict-Feedforward Nonlinear Systems”, *IEEE Transactions on Automatic Control*, 55(2), 2010, 287-303.

[19] Lozano, R., P. Castillo, P. Garcia and A. Dzul, “Robust Prediction-Based Control for Unstable Delay Systems: Application to the Yaw Control of a Mini-Helicopter”, *Automatica*, 40(4), 2004, 603-612.

[20] Nešic, D., and L. Grüne, “Lyapunov-based continuous-time nonlinear controller redesign for sampled-data implementation”, *Automatica*, 41, 2005 1143–1156.

[21] Olaru, S., and, S.-I. Niculescu, “Predictive control for linear systems with delayed input subject to constraints”, *Proceedings of the 17th World Congress of the International Federation of Automatic Control*, Seoul, Korea, 2008.

[22] Simoes, C., H. Nijmeijer and J. Tsinias, “Nonsmooth Stabilizability and Feedback Linearization of Discrete-Time Nonlinear Systems”, *International Journal of Robust and Nonlinear Control*, 6(3), 1996, 171-188.

[23] Velasco-Villa, M., B. del-Muro-Cuellar and A. Alvarez-Aguirre, “Smith-Predictor Compensator for a Delayed Omnidirectional Mobile Robot”, *Proceedings of the 15th Mediterranean Conference on Control and Automation*, Athens, Greece, 2007.